\documentclass{elsart}
\usepackage{graphics}
\begin{document}
\runauthor{Gh. Adam, S. Adam, and N.M. Plakida}
\begin{frontmatter}
\title{Reliability Conditions in Quadrature Algorithms}
\author[plak,adam]{Gh. Adam\thanksref{gemail}}
\author[plak,adam]{S. Adam\thanksref{semail}}
\author[plak]{N.M. Plakida\thanksref{pemail}}
\thanks[gemail]{Corresponding author;~e-mail: {\sf adamg@theory.nipne.ro}}
\thanks[semail]{~e-mail: {\sf adams@theory.nipne.ro}}
\thanks[pemail]{~e-mail: {\sf plakida@thsun1.jinr.ru}}
\address[plak]{Bogolubov Laboratory of Theoretical Physics,
Joint  Institute  for Nuclear Research, 141980 Dubna, Russia}
\address[adam]{Department of Theoretical Physics,
Institute of Physics and Nuclear Engineering,
P.O.~Box MG-6, 76900 Bucharest-M\u agurele, Romania}
\begin{abstract}
The detection of insufficiently resolved or ill-conditioned
integrand structures is critical for the reliability assessment
of the quadrature rule outputs. We discuss a method of analysis of the
{\sl profile of the integrand at the quadrature knots\/} which allows
inferences approaching the theoretical 100\% rate of success,
under error estimate sharpening.
The proposed procedure is of the highest interest for
the solution of parametric integrals arising in complex physical
models.
%
%
\end{abstract}
\end{frontmatter}
%
\section{Introduction}
\label{sec:intro}

A large number of physical models currently under study are characterized by
two combined features. First, the observables are obtained as integrals which
cannot be solved analytically. Second, the models describe physical systems
involving one or more specific parameters the variation of which results in
critical modification of the system behaviour.
  As a consequence, deep understanding of the predictions of the models needs
the exploration of the values of the observables over a large range of the
variable parameters.

As usual, to solve the occurring parametric integrals, recourse is made to
existing library codes of automatic adaptive quadrature which may fail
badly without providing any hint about such possibilities.
  We are directly aware of three such frustrating experiences.
  The first one concerns the two-band singlet-hole Hubbard model of cuprate
superconductors
\cite{Plak95}--\cite{Plak01},
which involves integrals over ranges of the first Brillouin zone. The
variation of the parameter of the model (the hole or electron doping in the
high-$T\sb{c}$ superconductor) results in substantial modification of
the behavior of the involved functions over the Brillouin zone.
The exploration of the predictions of the physical model with the doping is
fundamental for the validation of the proposed mechanism as responsible for
the superconducting pairing in cuprates. However, the reliability of the
outputs was found to be exceedingly low to allow sound inferences based on
the bare numerical outputs.
  A similar problem arises in the alternative
${{\mbox{U(1)}} \times {\mbox{SU(2)}}}$ gauge theory model of underdoped
cuprate superconductors
\cite{YuLu98}.
The meaningful physical solution derived under simplifying assumptions in
\cite{YuLu00}
could not be recovered from outputs generated by the available automatic
adaptive quadrature codes.
  The numerical exploration of a model of nuclear fission
\cite{Aurel98}
could not be achieved by means of library quadrature codes either.

These circumstances come from the fact that the existing algorithms for the
numerical integration of real valued functions (see, e.g.,
\cite{K&U98}
for details on the available algorithms and a recent review of numerical
quadrature) are tailored for specific classes of integrands, with limited
possibilities to solve simultaneously families of integrals falling in
different classes.

We may therefore assume that a study able to increase the reliability of the
automatic adaptive quadrature algorithms for solving parametric integrals in
connection with the exploration of physical models is of interest for a great
many users.
  Within the automatic adaptive quadrature, the approximate value $Q$ of a
given integral as well as its associated error estimate~$E$ are obtained as
sums of {\sl local couples\/} \{$q, e$\} of estimates over subranges.

The general picture offered by the numerical evidence on the solution of
parametric integrals points towards the existence of a limited range of
parameter values where the local quadrature sum $q$ provides accurate
solution of the integral of interest, whereas for other parameter values
the quadrature sum $q$ is inaccurate.
  Over the range of accurate $q$ outputs, the existing quadrature error
estimators provide outputs~$e$ which, in most cases, {\sl grossly
overestimate\/} the actual quadrature error, whence the need of
supplementary range subdivisions and over computing to meet the input
precision requests.
  However, over the range of inaccurate $q$ outputs, the heuristics
implemented in the local error estimators may result in spurious outputs
quoted as reliable, hence the impossibility to detect such cases by means
of the existing library codes.

In the present paper we discuss a generalization of the approach proposed in~%
\cite{AA01}
intended to reconcile these two contradictory aspects.
  The cornerstone of such an analysis is the derivation of reliability
criteria for the validation of the local error estimate $e$ associated to a
local quadrature sum $q$ based on the study of the {\sl profile of the
integrand at the set of quadrature knots} entering the expression of $q$.

The basic idea is that an unreliable estimate of $e$ might originate either
in the {\sl insufficient resolution\/} of the integrand profile, or in the
presence of {\sl difficult isolated points\/} (integrable singularities,
turning points, jumps) which result in slow convergence.
  The occurrence of each kind of difficulty can be evidenced by means of
specific consistency criteria asking for the fulfillment of requirements
following from quite general considerations: the very definition of the
Riemann integral, the fundamental properties of the basis polynomials
which span the approximating linear space where the interpolatory polynomial
of the quadrature rule is defined, the properties of the continuous
functions at or near their extremal points, and the smoothness properties
of the continuous functions inside their monotonicity subranges.

If the integrand is well-conditioned but its profile is
insufficiently resolved at the current set of quadrature knots,
repeated subdivision of the integration range eventually results in the
fulfillment of all the reliability constraints.
  A genuine difficult integrand point, however, recurs under repeated subrange
subdivisions.
  Therefore, repeated analysis of the integrand profile under subrange
subdivision ultimately results in the {\sl diagnostics stability under
iteration}. This is the point where the general control routine of an
automatic quadrature rule can take safe decisions concerning the best way
to continue the solution refinement or to decide that the integral was
solved within the input accuracy specifications.

The paper starts with basic definitions and notations (section~\ref{sec:defs}).
   In section~\ref{sec:stabl}, the main features of the validation procedure
of a computed local couple $\{q, e \}$ are discussed.
   Criteria for the identification of ill-conditioning features within an
integrand profile are summarized in section~\ref{sec:wellco}.
   Their practical importance is assessed in the section~\ref{sec:applics}
based on numerical evidence obtained from the solution of case study
integrals by Gauss-Kronrod 10--21 quadrature rules
\cite{QUADPACK}
with improved error estimate
\cite{AA01}.
   Concluding comments are given in section~\ref{sec:concl}.
\section{Definitions and notations}
\label{sec:defs}

Let $I$ denote the actual value of the integral to be solved numerically,
\begin{equation}
   I \equiv I[f] = \int\sb{a}\sp{b} g(x) f(x) dx \ ,\
     - \infty < a < b < \infty .
  \label{eq:Iref}
\end{equation}
Here, the {\sl weight function} $g(x)$ is an analytically integrable function
which absorbs a {\sl difficult part\/} of the integrand (e.g., an oscillatory
or a singular factor). In the absence of such factors, $g(x) = 1$.
The integrand function $f(x)$ is assumed to be {\sl continuous almost
everywhere\/} on $[a, b]$, such that~(\ref{eq:Iref}) exists and is finite.

A {\sl local quadrature rule\/} produces as solution of~(\ref{eq:Iref})
a couple \{$q,~e$\}, where the {\sl quadrature sum} $q$ yields an approximate
value of the integral $I$, while the {\sl local error estimate} $e > 0$
provides information on the accuracy of $q$. If $e > |e\sb{Q}|$, where
\begin{equation}
   e\sb{Q} = I - q 
  \label{eq:exct}
\end{equation}
is the actual error associated to $q$, then the couple
\{$q,~e$\} is {\sl reliable}, otherwise it is {\sl unreliable} and
the numerical solution fails.

A $(2n + 1)$-knot interpolatory quadrature sum $q\sb{2n}$ is obtained as
the analytical solution of the integral~(\ref{eq:Iref}) with the
integrand~$f(x)$ replaced by an interpolatory polynomial of the $2n$-th
degree,
\begin{equation}
   P\sb{2n}(x) = \sum\sb{k=0}\sp{2n} \alpha\sb{k}~p\sb{k}(x),
  \label{eq:Pref}
\end{equation}
where \{$p\sb{k}(x)$\} is the set of polynomials of degree at most
$2n$ spanning the approximating space of $P\sb{2n}(x)$. The coefficients~%
$\alpha\sb{k}$ are obtained from the set of conditions of interpolation
\begin{equation}
   P\sb{2n}(x\sb{i}) = f(x\sb{i}),
  \label{eq:icnd}
\end{equation}
at a set of $2n+1$ abscissas (called {\sl quadrature knots}) inside $[a, b]$,
\begin{equation}
   a\leq x\sb{0} < x\sb{1} < \cdots < x\sb{2n} \leq b\, .
  \label{eq:qkgen}
\end{equation}

In the particular case of the {\sl symmetric\/} ($2n+1$)-knot quadrature
sums, the interpolation abscissas inside $[a, b]$ are given by
\begin{equation}
   x\sb{i} = c + h y\sb{i};\; c = (b+a)/2;\; h = (b-a)/2;\;
         i = -n, -n+1, \cdots, n,
  \label{eq:qknots}
\end{equation}
where the reduced quadrature knots $y\sb{i}$ are defined on $[-1, 1]$, such
that
$
   0 = y\sb{0} < y\sb{1} < y\sb{2} < \cdots < y\sb{n} \leq 1,
$
while $y\sb{-i} = - y\sb{i},\; i = 1, \cdots, n$.

The local quadrature sum $q\sb{2n}$ is then expressed as a linear
combination of the integrand values at the quadrature knots,
\begin{equation}
   q\sb{2n} \equiv Q\sb{2n}[f] = \sum\sb{i = -n}\sp{n} w\sb{i} f(x\sb{i})\, ,
  \label{eq:rK}
\end{equation}
with the quadrature weights showing the symmetry property $w\sb{-i} = w\sb{i}$.

The information provided by the $2n+1$ integrand values at the quadrature
knots, $\{ f(x\sb{i}) | i = -n, \cdots, n \}$, is insufficient for the
derivation of an expression for the error estimate $e\sb{2n}$ associated
to $q\sb{2n}$.

Kronrod~%
\cite{Kronrod65}
derived an error estimate, called in what follows {\sl genuine
Gauss-Kronrod\/} ({\tt ggk}) error estimate, from an upper bound of
\begin{equation}
   e\sb{ggk} = \vert q\sb{2n} - q\sb{n} \vert \, ,
  \label{eq:eggk}
\end{equation}
where $q\sb{2n}$ is the quadrature sum~(\ref{eq:rK}), while~$q\sb{n}$ is a
lower degree quadrature sum derived over the subset of~(\ref{eq:qknots}),
\begin{equation}
   x\sb{-n+\gamma} < x\sb{-n+\gamma + 2} < \cdots < x\sb{n-\gamma - 2} <
   x\sb{n-\gamma}\, .
  \label{eq:qcoarse}
\end{equation}
Here, $\gamma = 1$ for an open quadrature sum (typically, the Gauss-Kronrod
({\tt GK}) quadrature where the spanning basis~$\{p\sb{k}(x)\}$
in~(\ref{eq:icnd}) is given by Legendre polynomials and their orthogonal
Kronrod extensions), while $\gamma = 0$ for a closed quadrature sum
(typically, the Clenshaw-Curtis ({\tt CC}) quadrature where the spanning
basis~$\{p\sb{k}(x)\}$ in~(\ref{eq:icnd}) is given by Chebyshev polynomials).

In what follows, the set of quadrature knots~(\ref{eq:qknots}) is referred to 
as the {\sl fine\/} discretization of the integration domain $[a, b]$, while
the sparser set of quadrature knots~(\ref{eq:qcoarse}) as the {\sl coarse\/}
discretization of $[a, b]$. The integrand values at these knots define its
{\sl fine\/} and {\sl coarse samplings\/} respectively.

In the {\tt QUADPACK} package
\cite{QUADPACK},
which has been incorporated in most major program libraries, while a {\tt ggk}
error estimate was implemented for the {\tt CC} quadrature, the {\tt GK} error
estimate was reformulated as follows.
   Let ${\bar {f}}$ denote the computed value of the average of $f(x)$ over
$[a, b]$ at the knots~(\ref{eq:qknots}),
\begin{equation}
   \bar {f} = q\sb{2n}/(b - a)\, ,
  \label{eq:barf}
\end{equation}
and let $\Delta = Q\sb{2n}\Bigl[ \vert f - \bar {f} \vert \Bigr]$ denote
the computed value of
$\int\sb{a}\sp{b} \vert f(x) - \bar {f} \vert dx$,
which measures the area covered by the deviations of $f(x)$ around~$\bar {f}$.

The local {\sl {\tt QUADPACK} error estimate\/} ({\tt qdp}) is then given by
\begin{equation}
   e\sb{qdp} = \Delta \times \min \{(200 e\sb{ggk}/\Delta)\sp{3/2},\; 1\}.
  \label{eq:eqdp}
\end{equation}

The values (\ref{eq:eggk}) and (\ref{eq:eqdp}) are taken for error estimates
provided they do not fall below the attainable accuracy limit imposed by the
relative machine precision. The latter threshold is defined as the product
\begin{equation}
   e\sb{roff} = \tau\sb{0} \epsilon\sb{0} Q\sb{2n}\Bigl[ \vert f \vert \Bigr].
  \label{eq:erof}
\end{equation}
Here $\tau\sb{0}$ is an empirical multiplicative factor (following
{\tt QUADPACK}, we have chosen $\tau\sb{0} = 50$) and $\epsilon\sb{0}$
denotes the relative machine accuracy.

For the case study integrals considered below, the value $I$ of~(\ref{eq:Iref})
is computed from the existing analytical expressions, such that the
{\sl exact\/} error $e\sb{Q}$~(\ref{eq:exct}) of the quadrature sum~$q\sb{2n}$
can be defined.

In the graphical representation of the quadrature errors, the {\sl moduli
of the relative errors\/} (simply called relative errors in the sequel) are
useful,
\begin{equation}
   \varepsilon \sb{\alpha} = \vert e\sb{\alpha} / I \vert \,,
     \quad \alpha \in \{ 2n, Q \}.
  \label{eq:erel}
\end{equation}

The derivation of the local error estimates within a subroutine
which implements a quadrature rule uses information inferred from
the {\sl estimated relative errors},
\begin{equation}
   \rho \sb{\alpha} = \vert e\sb{\alpha} / q\sb{2n} \vert \, ,
     \quad \alpha \in \{ ggk, qdp, 2n \}.
  \label{eq:cerel}
\end{equation}
\section{Stability of the diagnostics under subrange subdivision}
\label{sec:stabl}

Using~(\ref{eq:eggk}),~(\ref{eq:eqdp}) and~(\ref{eq:erof}), we get the local
error estimate
\cite{AA01}
\begin{equation}
  e\sb{2n} = \max\left[
                     e\sb{roff},\,
                     \min (e\sb{ggk}, e\sb{qdp})
                 \right]\, ,
  \label{eq:practic}
\end{equation}
the reliability of which is almost always subject to doubt, except for the
case when the lower degree quadrature sum $q\sb{n}$ is sufficiently
accurate such that the accuracy of $q\sb{2n}$ itself reaches nine to ten
significant figures at least.
Such a condition can be confidently assumed to hold provided
\begin{equation}
  e\sb{2n} > e\sb{thr},\quad e\sb{thr} = 2\sp{-18} \simeq 0.38\times10\sp{-5}.
  \label{eq:rlb0}
\end{equation}
This empirically proposed threshold value is about two decimal figures
more conservative than the smallest values of the unreliable computed error
estimates over the evidence discussed in Sec.~\ref{sec:applics}.

If the opposite of~(\ref{eq:rlb0}) occurs, then a {\sl validation
procedure\/} is to be used to assess the reliability of the local couple
\{$q\sb{2n}, e\sb{2n}$\}. Thus, a self-validating quadrature rule returns,
besides the numerical output for $e\sb{2n}$, a flag having the value zero
in case of assumed reliable outputs and non-zero value if the output
is not validated.

The validation procedure proposed in this paper is based on the study of the
information contained in the {\sl profile of the integrand $f(x)$}
over $[a, b]$, defined as the set of integrand values at the quadrature
knots~(\ref{eq:qkgen}), completed with the endpoint values $f(a)$ and $f(b)$
in the case of open quadrature sums.
  Since the operation of subrange subdivision within automatic adaptive
quadrature always involves inner abscissas at existing quadrature knots,
the only price to be paid for the inclusion of the endpoint values in the
integrand profile is the direct access of the general control routine to
such data. This goal is achieved provided the generation of the integrand
sampling at the quadrature knots~(\ref{eq:qkgen}) is done within a
subroutine which is distinct from that implementing the quadrature rule
and is directly subordinated to the general control routine.

The study of the integrand profile starts with the definition of its
{\sl monotonicity subranges},
$[x\sb{i\sb{j-1}}, x\sb{i\sb{j}}]$, over $[a, b]$, where
\begin{equation}
   a = x\sb{i\sb{0}} < x\sb{i\sb{1}} < x\sb{i\sb{2}} < \cdots < x\sb{i\sb{m}}
   < x\sb{i\sb{m+1}} = b\, ,
  \label{eq:mntsb}
\end{equation}
denote the abscissas of the extremal points of $f(x)$ within the sampling.

In terms of the answer concerning the number of monotonicity subranges,
several specific reliability criteria are checked and the number $\lambda$
of the infringements of these criteria is counted. There are three
critical values of the pointer $\lambda$ in terms of which a decision is
taken:
\begin{itemize}
  \item
      $\lambda = 0$: probably the investigated integrand profile comes
      from a well-conditioned integrand, hence the output~$q\sb{2n}$,
      Eq.~(\ref{eq:rK}) is reliable, while the quadrature error
      estimate~(\ref{eq:practic}) is overestimated.
  \item
      $\lambda = 1$ or $\lambda = 2$: there is a high probability that a
      difficult isolated point is present which implies slow convergence
      of the quadrature sums.
  \item
      $\lambda \geq 3$: the insufficient resolution of the integrand
      profile at the involved quadrature knots is manifest. The output
      is useless and further subrange subdivisions are compulsory.
\end{itemize}

The existence and finiteness of the Riemann integral~(\ref{eq:Iref})
guarantees that, after a {\sl finite\/} number of subrange subdivisions,
the discretization process will reach a {\sl stable\/} profile configuration
the refinement of which will result in unessential modifications only.

Under the occurrence of {\sl isolated\/} difficult points of the integrand,
the discretization process will resolve the profile over the well-conditioned
subranges within a finite number of subrange subdivisions, and then it
will mainly create a dense mesh around the difficult points.
  In this case, the automatic control subroutine will safely decide upon
the activation of a specific convergence acceleration algorithm, such
that a reliable numerical solution will be available in the end.

The achievement of the {\sl stability\/} of the diagnostics concerning the
conditioning properties of the integrand profiles over subranges, got after
a {\sl finite\/} number of subrange subdivisions, is the fundamental feature
which secures the efficiency of the procedure proposed in this investigation.

The occurrence of consistent with each other reliability diagnostics over
the current integration range and its subranges obtained by subrange
subdivision enables the general control routine to take safe decisions
concerning the activation of the implemented alternative algorithms.
\section{Well-conditioned integrand profiles}
\label{sec:wellco}

The consistency requirements satisfied by a well-conditioned integrand profile
are formulated mostly locally and they follow from quite general considerations
which are discussed in the next subsections.

Any infringement of the consistency criteria derived below is to be added
to the value of the ill-conditioning pointer $\lambda$.
\subsection{Insensitivity of the integral sums to discretization details}
\label{sec:Riemann}

The standard definition of the integral sums in a Riemann integral assumes
the fulfillment of the following two features:
\begin{itemize}
  \item[(i)]
       The norm of the discretization step defined over the integration domain
       tends to zero.
  \item[(ii)]
       The integral sum is insensitive to the the addition or removal of a
       single discretization abscissa within the defined partition.
\end{itemize}

In the quadrature algorithms, the norm of the discretization~(\ref{eq:qknots})
is far from being close to zero. The quadrature knots are not distributed
evenly either.
   For the {\tt GK} and {\tt CC} quadrature rules mentioned above, the
fundamental range $[-1, 1]$ consists of a sparser knot region centered
around the origin and two denser knot regions located toward the
range ends. The number of abscissas entering the integrand profile
associated to a ($2 n + 1$)-knot open quadrature rule equals $2 n + 3$,
while the corresponding number for a closed quadrature rule is $2 n + 1$.
Therefore, for both kinds of quadrature rules, a particular inner reduced
knot $y\sb{i}$ lies in the dense knot region provided the lengths of its
two adjacent subranges are smaller than the threshold quantity for a
uniform distribution, $d\sb{av} = 2/(2 n + 3)$.

An immediate consequence of the feature (i) is the property that the denser
discretization regions of a {\sl smooth\/} integrand $f(x)$ secure better
accuracy of their contributions to the quadrature sums than the sparser ones.
We reformulate this observation as follows: the generation of the fine
discretization~(\ref{eq:qknots}) from the coarse discretization~%
(\ref{eq:qcoarse}) is expected to result in non-essential modifications of
the profile of $f(x)$ over the regions of dense knot discretization.

To characterize the extent to which a profile is modified by the addition of
new knots inside the region of dense knot discretization, let us consider that
$x\sb{0}$ is such a knot.
  If $x\sb{0}$ belongs to the set of extremal points~(\ref{eq:mntsb}) such
that the integrand value $f(x\sb{0})$ is {\sl isolated\/} from the integrand
values $f(x\sb{-1})$ and $f(x\sb{1})$ at the nearest neighbours $x\sb{-1}$
and $x\sb{1}$ by the median line $f = \bar {f}$, Eq.~(\ref{eq:barf}),
then the knot $x\sb{0}$ is said to be {\sl sensitive}. If {\sl both\/}
quantities $f(x\sb{-1})$ and $f(x\sb{1})$ stay on the same side with
$f(x\sb{0})$ with respect to the median line $f = \bar {f}$, then the
knot $x\sb{0}$ is said to be {\sl regular}. If the median line
$f = \bar {f}$ separates $f(x\sb{0})$ from {\sl only one\/} of the
values $f(x\sb{-1})$ or $f(x\sb{1})$, then the knot $x\sb{0}$ is said to be
{\sl gray}.

   We are now ready to formulate the first practical reliability criterion:
\begin{itemize}
  \item[(I)]
      {\tt Non-sensitivity of the extremal points:}\\ 
      {\sl The addition of supplementary quadrature knots to the coarse
      partition~(\ref{eq:qcoarse}) to reach the fine partition~%
      (\ref{eq:qknots}) does not result in supplementary gray or sensitive
      extrema of the profile of a well-conditioned integrand over the
      regions of dense knot discretization}.
\end{itemize}
\subsection{Features which stem from the basis polynomials}
\label{sec:basicp}

Since the equations~(\ref{eq:qknots}) perform the mapping of the original
interval $[a, b]$ onto the reduced interval $[-1, 1]$ over which the
orthogonal polynomials are usually defined, in this subsection we refer
to this reduced interval and use the notation $p\sb{k}(y)$ for the
basis polynomials. All the properties discussed below hold true over
arbitrary interval lengths, hence reference to the expression~(\ref{eq:Pref})
of the interpolatory polynomial spanned by the basis orthogonal polynomials
does not give rise to any confusion.

The existence and uniqueness of the interpolatory polynomial~(\ref{eq:Pref})
is secured provided the set of basis polynomials spanning~(\ref{eq:Pref})
define a Chebyshev system over $[a, b]$. Therefrom the following properties
hold true:
\begin{itemize}
  \item[(iii)]
       $p\sb{0}(y) = const$.
  \item[(iv)]
       The set of the successive extremal values of a polynomial~$p\sb{k}(y)$
       of degree $k > 1$ defines an {\sl alternating sequence\/}
       over~$[-1, 1]$.
  \item[(v)]
       The zeros of the polynomials $p\sb{k}(y)$ and $p\sb{k+1}(y)$
       {\sl are interlaced\/} inside the open range $(-1, 1)$.
\end{itemize}

The average value ${\bar {f}}$, Eq.~(\ref{eq:barf}), of the
integrand~$f(x)$, which defines its zeroth order moment over the
sampling~(\ref{eq:qkgen}) and is related to the coefficient of
$p\sb{0}(y)$ within a basis set of orthogonal polynomials, serves
as reference value with respect to which the oscillations of the
integrand profile are counted.
The intersections of the integrand profile with the line
$f = {\bar {f}}$ define the {\sl zeros\/} of the integrand profile.

The alternation property (iv) results in the important consequence that
the deviations of the successive extremal values of a well-conditioned
integrand profile from~${\bar {f}}$ define an {\sl alternating\/}
sequence with {\sl comparable amplitudes\/} at the adjacent extremal
knots~(\ref{eq:mntsb}). This property can be detailed for practical
purposes into two well-conditioning alternation criteria:
\begin{itemize}
 \item[(IIa)]
  {\tt Type-1 alternation criterion:}\\
  -- {\sl Each inner monotonicity subrange of a well-conditioned integrand
     profile\\
     \phantom{ -- }intersects the line $f = \bar{f}$.\\
  -- The two end point monotonicity subranges do not diverge
     from $f = \bar{f}$}.
 \item[(IIb)]
  {\tt Type-2 alternation criterion:}\\
     {\sl Each inner gray extremal point which satisfies the type--1
     alternation criterion is to stay sufficiently far from the line
     $f = \bar{f}$}.
\end{itemize}

The test for the occurrence of an infringement of the type--1 alternation
criterion is obvious. As it concerns the the latter criterion, two
infringements are to be simultaneously tested:
\begin{itemize}
 \item
     The distance from $f(x\sb{0})$ to $\bar{f}$ is to be smaller than those
     of its nearest neighbouring extrema.
 \item
     Let $a\sb{0}$, $a\sb{l}$ and $a\sb{r}$ denote the areas
     surrounded by $f = \bar{f}$ and the integrand profile around $x\sb{0}$
     and its nearest neighbours in the set~(\ref{eq:mntsb}). Then
     \begin{equation}
        |a\sb{0}| < t\sb{1}|a\sb{l} + a\sb{r}|,\quad t\sb{1} = 10,
      \label{eq:critIIb}
     \end{equation}
     where the value of $t\sb{1}$ was chosen such as to point to a
     discrepancy exceeding an order of magnitude.
     The computation of the three local areas is done by compound trapeze rule
     which is robust and sufficiently accurate for the involved comparison.
\end{itemize}

Corroboration of the interlacing property (v) with the non-sensitivity
criterion (I) results in a criterion for the distribution of the zeros
of the integrand profile:
\begin{itemize}
  \item[(III)]
      {\tt Non-sensitivity of the zeros:}\\
      {\sl Over the dense knot regions, the numbers of zeros of the fine
      and coarse profiles of a well-conditioned integrand are the same.}
\end{itemize}
\subsection{Integrand variations around its isolated extremal points}
\label{sec:xisola}

The lateral first order derivatives of a smooth first order differentiable
function vanish at an extremal point, while the curvature of a second order
differentiable function (which is given by the second order derivative)
keeps constant sign over a nonvanishing neighbourhood of the extremum.

Within the discrete mesh defined by the quadrature knots, inquiries about
these properties can be made only at integrand profile approximations of
{\sl isolated\/} extremal points
of the integrand. If $x\sb{0}$ is such a point, then a {\sl sufficiently
large\/} neighbourhood $\{ \xi\sb{l}, \xi\sb{r}\}$ around $x\sb{0}$ can be
defined within which the evaluation of the quantities of interest is
expected to be weakly influenced by the presence of neighbouring extrema.

Let us assume that an isolated extremal point of a well-conditioned integrand
was identified within a sufficiently well resolved integrand profile.
The following consistency criteria establish well-conditioned behaviours
of the data:
\begin{itemize}
  \item[(IV)]
      {\tt First lateral derivative criterion:}\\
       {\sl The approximation of the lateral first order derivatives at an
       isolated extremum of the profile using fine sampling data is closer
       to zero as compared to the value estimated from data defined over a
       coarse sampling with respect to the extremum location}.
  \item[(V)]
      {\tt Curvature sign constancy criterion:}\\
       {\sl The sign of the second order derivative computed from fine sampling
       data centered at the extremum is the same as that of the value
       estimated from data involving the reference extremum as a lateral
       point to the left/right}.
\end{itemize}

We shall illustrate the quantitative implementation of these criteria for
a reference extremum $x\sb{0}$ which is said to be {\sl isolated to the right}.
That is, the neighbourhood $\{ \xi\sb{l}, \xi\sb{r}\}$ contains inside it
the set of abscissas $\{ x\sb{-1}, x\sb{0}, x\sb{1}, x\sb{2} \}$
at which the integrand function takes respectively the values
$\{ f\sb{-1}, f\sb{0}, f\sb{1}, f\sb{2} \}$.

To estimate the approximation of the first order right lateral derivative,
we define the interpolatory polynomial of the third degree which fits these
four data. This yields the following result:
\begin{equation}
   f\sp{\prime}\sb{r, fine}(x\sb{0}) = d\sp{(1)}\sb{1,0} -
     \frac{h\sb{1,0}}{h\sb{2,-1}}\Big[ h\sb{0,-1} d\sp{(2)}\sb{2,1} +
            h\sb{2,0} d\sp{(2)}\sb{1,-1} \Big]\, .
  \label{eq:d1rf}
\end{equation}
Here, $h\sb{i,j} = x\sb{i} - x\sb{j}$,
      $d\sp{(1)}\sb{i,j} = (f\sb{i} - f\sb{j})/h\sb{i,j}$ denote the
first order divided differences at $x\sb{i}$ and $x\sb{j}$,
while $d\sp{(2)}\sb{2,1} = \left( d\sp{(1)}\sb{2,0} - d\sp{(1)}\sb{1,0}
\right) / h\sb{2,1}$ and
      $d\sp{(2)}\sb{1,-1} = \left( d\sp{(1)}\sb{1,0} - d\sp{(1)}\sb{0,-1}
\right) / h\sb{1,-1}$ denote specific second order divided differences.

On the other hand, the coarse sampling around $x\sb{0}$ yields:
\begin{equation}
   f\sp{\prime}\sb{r, coarse}(x\sb{0}) = d\sp{(1)}\sb{2,0}\, .
  \label{eq:d1rc}
\end{equation}
The criterion (IV) then simply states that the approximations~(\ref{eq:d1rf})
and~(\ref{eq:d1rc}) should satisfy
$\vert f\sp{\prime}\sb{r, fine}(x\sb{0}) \vert <
 \vert f\sp{\prime}\sb{r, coarse}(x\sb{0}) \vert$.

Over the same set of data, the criterion (V) requirement of constancy of the
sign of the second order derivative results in the condition
\begin{equation}
  \left( d\sp{(1)}\sb{2,0} - d\sp{(1)}\sb{1,0} \right)
  \left( d\sp{(1)}\sb{1,0} - d\sp{(1)}\sb{0,-1} \right) > 0\, .
  \label{eq:d2r}
\end{equation}

For the extremal point $x\sb{0}$ isolated to the left, similar conditions are
derived from the data set $\{ f\sb{-2}, f\sb{-1}, f\sb{0}, f\sb{1} \}$
obtained at the abscissas $\{ x\sb{-2}, x\sb{-1}, x\sb{0}, x\sb{1} \}$.
\subsection{Well-conditioning inside monotonicity subranges}
\label{sec:monowc}

Inside any monotonicity subrange of a smooth first order differentiable
function~$f(x)$, the first order derivative varies {\sl smoothly\/} from point
to point.

Within numerical quadrature, the fulfillment of this property for an
integrand sampling can be checked by making use of first order
divided differences.
If the integrand profile is monotonic over $[a, b]$, or monotonicity
subranges can be defined which extend over three successive knots at
least, then a smoothly varying profile will by characterized by the
{\sl absence of jumps:}
\begin{itemize}
  \item[(VI)]
      {\tt Absence of jumps inside monotonicity subranges:}\\
       {\sl Inside a monotonicity range, the ratio of two successive first
       order divided differences cannot exceed a relative smoothness threshold}.
\end{itemize}
If one of the knots involved in the divided differences is an extremal point,
then this smoothness condition is to be checked only one-directionally,
skipping the case of vanishingly small divided difference at the extremal
point.

For knots far from inflection points, a threshold value $t\sb{jmp} = 10$,
corresponding to the agreement of the successive divided differences within
an order of magnitude, is appropriate. In the neighbourhood of inflection
points characterized by a maximum of the first order derivative, this value
is to be halved to detect ill-conditioned behaviour, while in the
neighbourhood of inflection points characterized by a minimum of the first
order derivative, five times larger threshold value is appropriate.
\section{Numerical results}
\label{sec:applics}

The significance of the conditioning criteria discussed in the previous
section is intuitive and straightforward. In addition to the case specified
by the condition~(\ref{eq:rlb0}), a second case when the reliability analysis
can be skipped is that of a {\sl monotonic profile\/} characterized by an
error estimate
\begin{equation}
  e\sb{2n} > 0.5 |q\sb{2n}|\, .
 \label{eq:0.5}
\end{equation}
Then the computed quadrature sum is highly inaccurate, such that an error flag
can be directly assigned.

The diagnostics of the reliability criteria (IIb), (IV), (V), and (VI) depend
on specific adjustable parameters. If the quantitative thresholds entering
these criteria are decreased, the diagnostics becomes less permissive, with the
consequence that the reliability range shrinks and the number of wrong
diagnostics is decreased. The opposite occurs under the increase of the
quantitative thresholds. The numerical data reported in this section
show that, when corroborated with the requirement of the stability of the
diagnostics formulated in Sec.~\ref{sec:stabl}, the formulation of the
reliability criteria in Sec.~\ref{sec:wellco} is able to eliminate practically
all the spurious outputs occurring in an automatic adaptive quadrature
algorithm.

To illustrate the present analysis, a comparison is done of three codes using
Gauss-Kronrod 10--21 ({\tt GK 10-21}) quadrature rules:
($\alpha$) the {\tt QUADPACK} code~%
\cite{QUADPACK},
($\beta$) the self-validating code of~%
\cite{AA01},
and ($\gamma$) the code using the present reliability analysis.

Each code solved the parametric families of elementary integrals considered
in ref.~%
\cite{AA01}.

The first is the family of integrals over $[0, 1]$ of the terms of the
fundamental power series, $x\sp{n}$,
\begin{equation}
  \int\sb{0}\sp{1} x\sp{n}\,dx~= \frac{1}{n+1}\,, \quad n = 0, 1,\cdots , 1023.
 \label{eq:pow}
\end{equation}
\begin{figure}[ht]
\resizebox{\textwidth}{!}
{\includegraphics{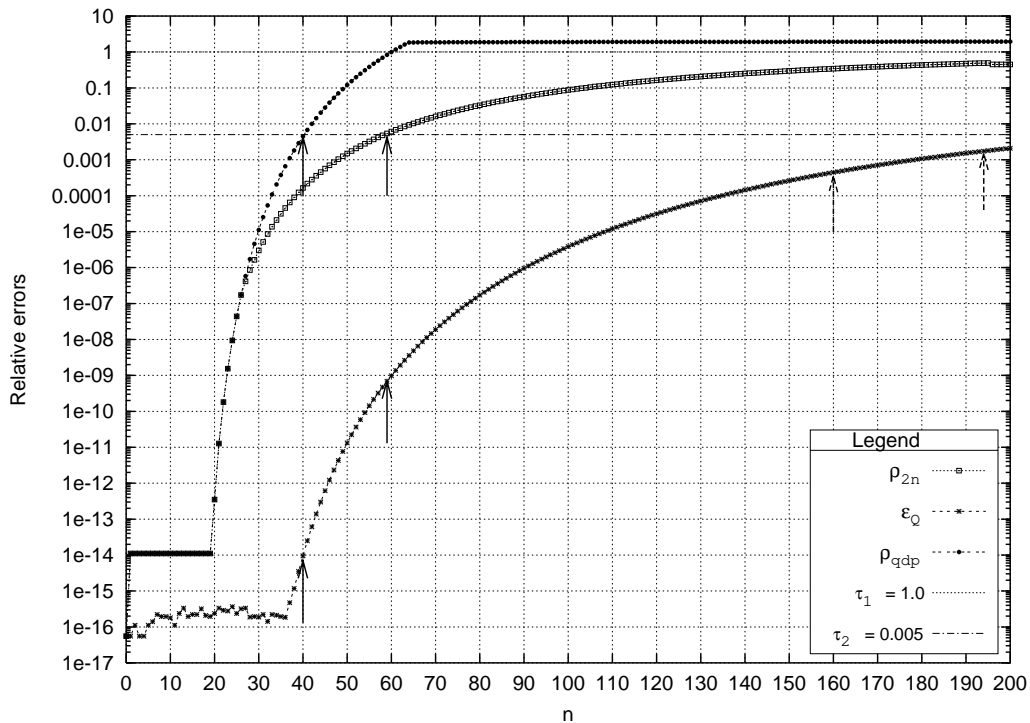}}
\caption{Relative errors $\rho\sb{qdp}$ and $\rho\sb{2n}$, Eq.~(\ref{eq:cerel}),
         $\varepsilon\sb{Q}$, Eq.~(\ref{eq:erel}), of the {\tt GK 10-21}
         outputs for the family of integrals~(\ref{eq:pow}) at exponents
         $n \leq 200$.
	 The upper leftmost solid line arrow points to the accuracy basin
         extension established by the {\tt QUADPACK} code using the error
	 estimate~(\ref{eq:eqdp}). The solid line arrow on the same vertical
	 points to the upper accuracy of the quadrature sum $q\sb{2n}$
	 retained as reliable by the {\tt QUADPACK} code.
	 The next pair of solid arrows show the result of the analysis done
	 in ref.~%
	 \cite{AA01}.
	 The left interrupted line arrow represents the extension of the
	 reliability basin established by the present analysis, while the
	 right one shows the exponent threshold above which the criterion~%
	 (\ref{eq:0.5}) supersedes the need of reliability analysis.
        }
\label{fig:fr21pow}
\end{figure}
The integrands are monotonic, inflection points are absent over the integration
range. Fig.~\ref{fig:fr21pow} illustrates the behaviour of the error estimates
with the power $n$ running over the range \{$0, 200$\}.
The results obtained for this family of integrals can be summarized as follows:
\begin{itemize}
 \item
   The {\tt QUADPACK} code infers an accuracy basin of the {\tt GK 10-21}
   code extending from $n = 0$ to $n\sb{max} = 40$, with the consequence that
   all the $q\sb{2n}$ outputs showing an actual accuracy lower than 14 decimal
   digits are thrown away. As shown in~%
   \cite{AA01},
   this early cut of the accuracy basin {\sl does not rule out the
   possibility of wrong error diagnostics\/} at asymptotically large $n$.
 \item
   The self-validating analysis of ref~%
   \cite{AA01}
   extends the accuracy basin of {\tt GK 10-21} up to $n\sb{max} = 59$,
   which corresponds to a correct identification
   of the outputs $q\sb{2n}$ as reliable up to accuracies of roughly nine
   significant digits. Above $n \ge 60$, {\sl all\/} the reliability
   diagnostics have been correct.
 \item
   The present analysis establishes an accuracy basin up to $n\sb{max} = 160$,
   which corresponds to outputs $q\sb{2n}$
   showing at least three significant decimal digits.
   At exponents $n \ge 195$, the criterion~(\ref{eq:0.5}) directly
   establishes the occurrence of unreliable $q\sb{2n}$ outputs without
   making recourse to the reliability analysis.
\end{itemize}

The second family solves integrals for a same integrand (which simulates a
centrifugal potential at large $x$) over ranges of variable length,
\begin{equation}
  \int\sb{0}\sp{b} \frac{1}{x\sp{2} + 1}\,dx~= \arctan (b)\,,
           \quad b = n,\quad n = 0, 1,\cdots , 10000.
 \label{eq:atg}
\end{equation}
\begin{figure}[h]
\resizebox{\textwidth}{!}
{\includegraphics{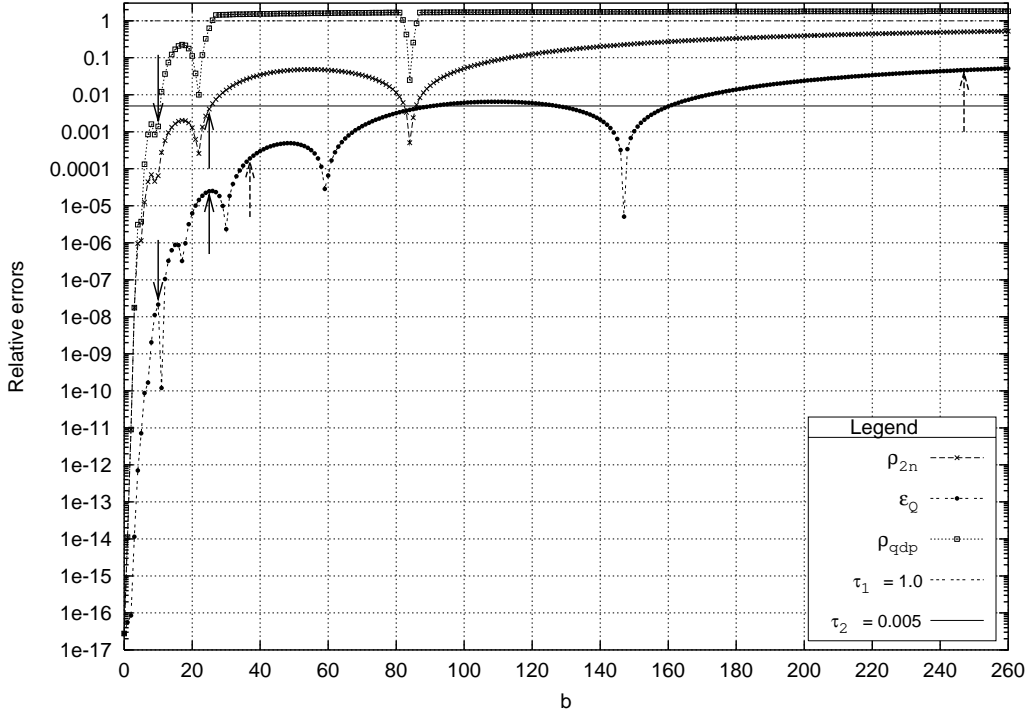}}
\caption{Same as fig.~\ref{fig:fr21pow} for the family of
         integrals~(\ref{eq:atg}), at upper integration ranges $b \leq 260$.
	 The arrows bear the same significance.
        }
\label{fig:fr21atg}
\end{figure}
The integrands are monotonic, an inflection point is present.
Fig.~\ref{fig:fr21atg} illustrates the behaviour of the error estimates
with the upper integration range $b = n$ for $n$ running over the range
\{$0, 260$\}.
In this figure the occurrence of cusps in the $\varepsilon \sb{Q}$ curve
points to the existence of fractional integration domain lengths at which the
quadrature sum $q\sb{2n}$ solves {\sl exactly\/} the integral~(\ref{eq:atg}),
such that the exact error changes sign.
The results obtained for this family of integrals can be summarized as follows:
\begin{itemize}
 \item
   The {\tt QUADPACK} code infers an accuracy basin of the {\tt GK 10-21}
   code extending up to $n\sb{max} = 10$, with the consequence that
   all the $q\sb{2n}$ outputs showing an actual accuracy lower than 
   seven decimal digits are thrown away. {\sl All the {\tt QUADPACK}
   reliability diagnostics above $b = n = 2460$ have been false}.
 \item
   The self-validating analysis of ref~%
   \cite{AA01}
   extends the accuracy basin of {\tt GK 10-21} up to $n\sb{max} = 25$,
   which corresponds to a correct identification
   of the outputs $q\sb{2n}$ as reliable up to accuracies of about six
   significant digits. At $n \ge 26$, {\sl all\/} the reliability
   diagnostics have been correct.
 \item
   The present analysis establishes an accuracy basin up to $n\sb{max} = 37$,
   which corresponds to outputs $q\sb{2n}$
   showing at least three significant decimal digits.
   At exponents $n \ge 247$, the criterion~(\ref{eq:0.5}) directly
   establishes the occurrence of unreliable $q\sb{2n}$ outputs without
   making recourse to the reliability analysis.
\end{itemize}
\begin{figure}[h]
\resizebox{\textwidth}{!}
{\includegraphics{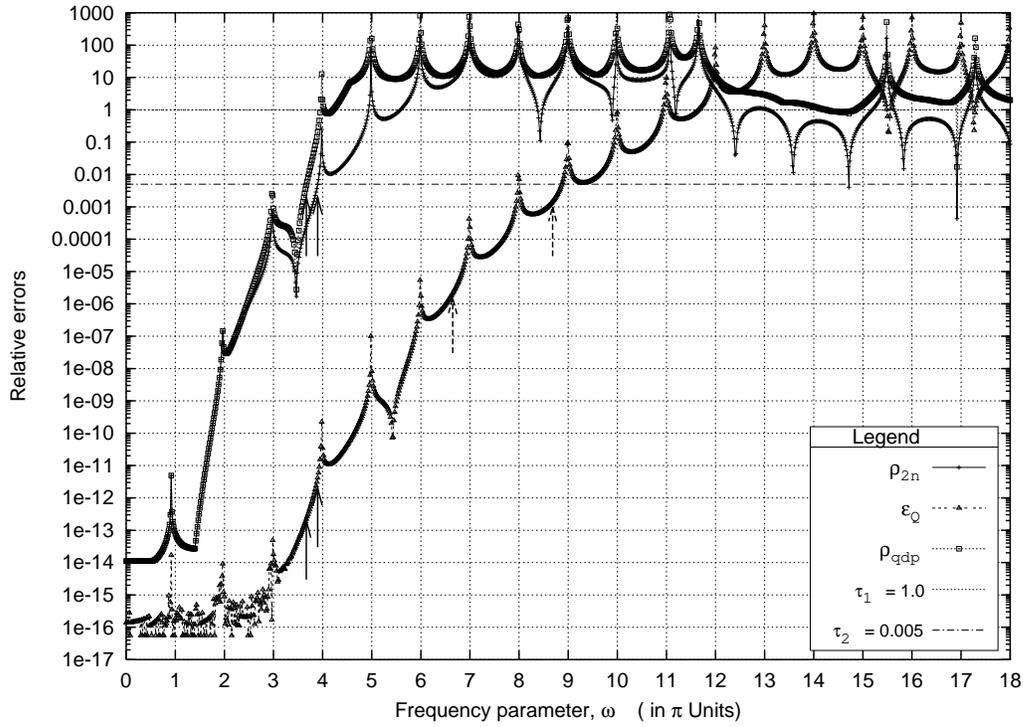}}
\caption{Outputs of the {\tt GK~10-21} quadrature rule for the family of
         integrals~(\ref{eq:oc10}) at $p = 1$.
         The significances of the solid line arrows are the same as in
	 fig.~\ref{fig:fr21pow}.
         The left interrupted line arrow represents the extension of the
	 reliability basin up to which the present analysis validates
	 {\sl all\/} the outputs~$q\sb{2n}$.
	 Inbetween the two interrupted line arrows, the diagnostics of the
	 present analysis is too conservative in about one third of the
	 solved cases.
        }
 \label{fig:oc10p1}
\end{figure}
\begin{figure}[ht]
\resizebox{\textwidth}{!}
{\includegraphics{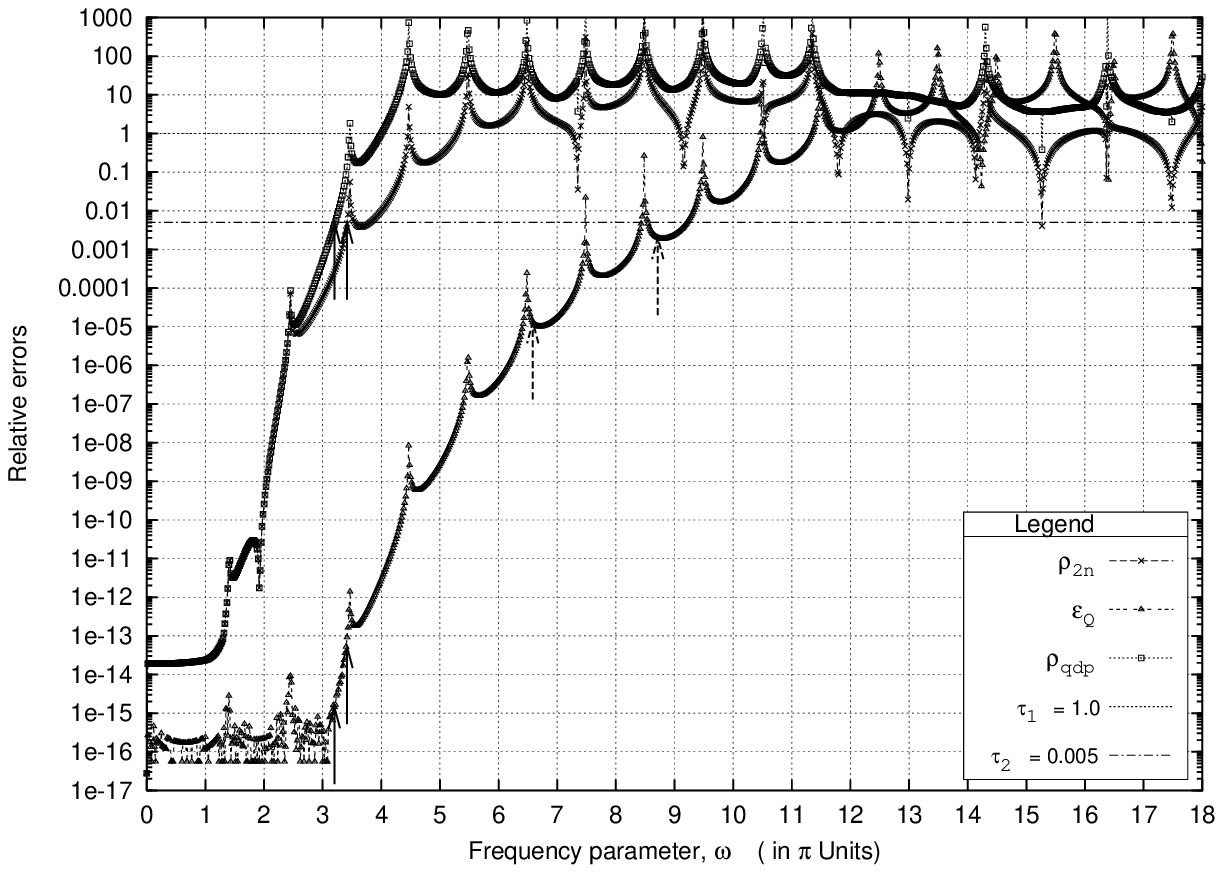}}
\caption{Same as fig.~\ref{fig:oc10p1} for the family of
         integrals~(\ref{eq:os10}) at $p = 1$.
        }
 \label{fig:os10p1}
\end{figure}

Next, we considered two pairs of families of integrals showing nonmonotonic
(oscillatory) behaviour, written in algebraically equivalent forms:
\begin{eqnarray}
  ({\bf C1})&&\int\sb{-1}\sp{1}e\sp{p(x-x\sb{0})}\cos(\omega x)\,dx~=
    \label{eq:oc10}\\
  ({\bf C2})&&\int\sb{0}\sp{1} 2e\sp{-px\sb{0}}\cosh(px)\cos(\omega x)\,dx~=
    \label{eq:oc11}\\
        && = 2 e\sp{-px\sb{0}}[p \sinh(p) \cos(\omega) +
                 \omega \cosh(p)\sin(\omega)] /(\omega\sp{2} + p\sp{2})\,;
    \label{eq:roc}
\end{eqnarray}
\begin{eqnarray}
  ({\bf S1})&&\int\sb{-1}\sp{1}e\sp{p(x-x\sb{0})}\sin(\omega x)\,dx~=
    \label{eq:os10}\\
  ({\bf S2})&&\int\sb{0}\sp{1} 2e\sp{-px\sb{0}}\sinh(px)\sin(\omega x)\,dx~=
    \label{eq:os11}\\
        && = 2 e\sp{-px\sb{0}}[p \cosh(p) \sin(\omega) -
                 \omega \sinh(p)\cos(\omega)] /(\omega\sp{2} + p\sp{2})\,.
    \label{eq:ros}
\end{eqnarray}

The parameter $\omega$ was chosen to run over the set of values
\begin{equation}
  \omega\sb{n} = n \pi /60,\quad n \in \{ 0, 6000\},
    \label{eq:omegan}
\end{equation}
while constant values $p = 1$  and $x\sb{0} = - 1$ have been chosen on the
ground that they are typical for the description of the behaviour of the
numerical results.
\begin{figure}[ht]
\resizebox{\textwidth}{!}
{\includegraphics{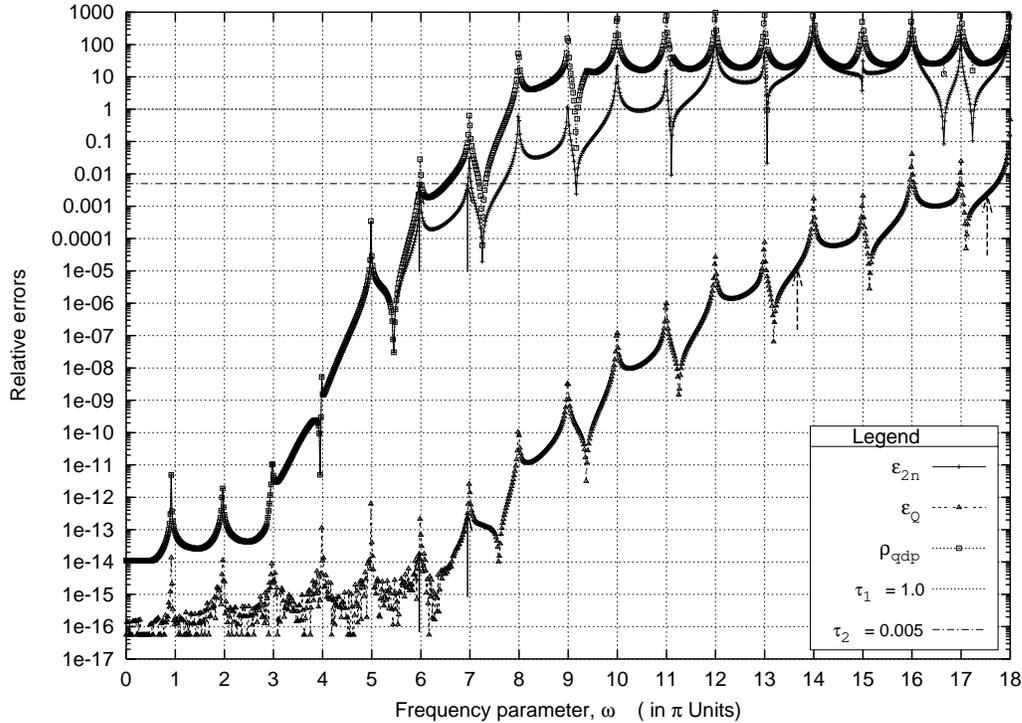}}
\caption{Same as fig.~\ref{fig:oc10p1} for the family of
         integrals~(\ref{eq:oc11}) at $p = 1$.
        }
 \label{fig:oc11p1}
\end{figure}
\begin{figure}[h]
\resizebox{\textwidth}{!}
{\includegraphics{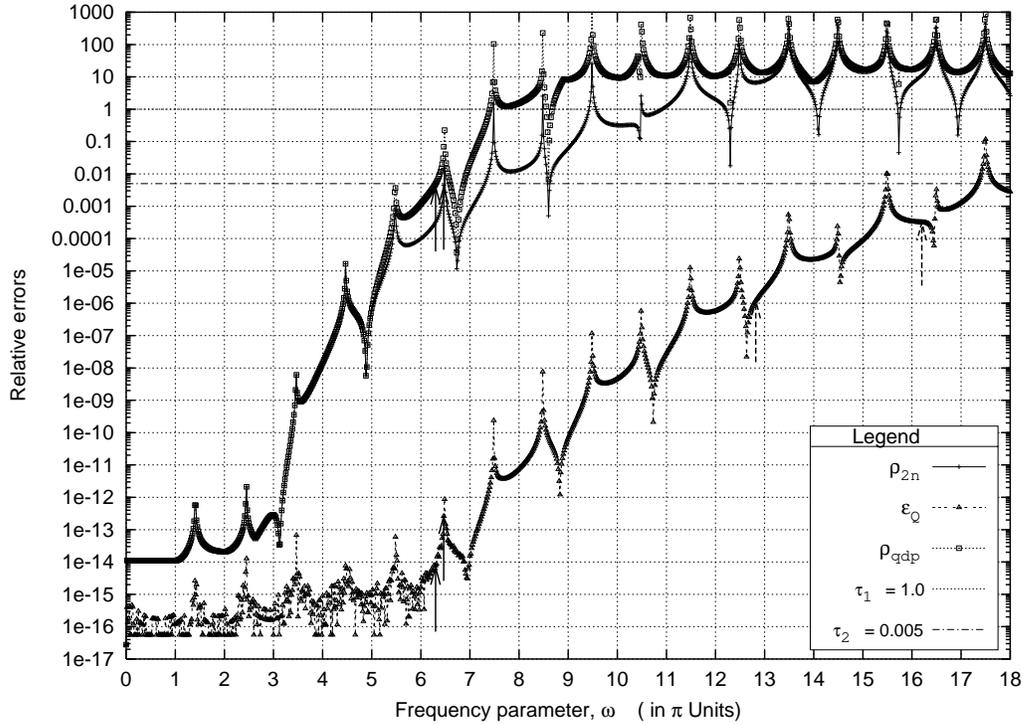}}
\caption{Same as fig.~\ref{fig:oc10p1} for the family of
         integrals~(\ref{eq:os11}) at $p = 1$.
        }
 \label{fig:os11p1}
\end{figure}

The analysis of the families of integrals~(\ref{eq:oc10}--\ref{eq:os11})
shows that the identification of a well-conditioned nonmonotonic integrand
profile needs testing the {\sl complete set\/} of consistency criteria
established in Sec.~\ref{sec:wellco}.
Therefore, the analysis is long. However, it is straightforward and
can be easily implemented in a computer program.

Figures~\ref{fig:oc10p1} to~\ref{fig:os11p1} show outputs
for the parameter $n$ running over
the range $n \in \{ 0, 1080\}$.
In these figures, two peculiarities of the $\varepsilon \sb{Q}$ curves are
apparent.
Similar to Fig.~\ref{fig:fr21atg}, the occurrence of cusps at
{\sl minima\/} in the $\varepsilon \sb{Q}$ curves point to the existence
of values of the parameter $\omega$ at which the given integrals are solved
{\sl exactly\/} by the quadrature sum $q\sb{2n}$, such that the exact error
changes sign.
The sharp maxima noticed in the $\varepsilon \sb{Q}$ curves occur at $\omega$
values which correspond to {\sl entire periods\/} of the oscillatory
factors over the integration range, such that important {\sl cancellation
by subtraction\/} effects occur which result in sensible worsening of the
numerical output.

A summary of the results obtained for the families of integrals~(\ref{eq:oc10}),
(\ref{eq:oc11}), (\ref{eq:os10}), and~(\ref{eq:os11}) is given in Table~1.

The {\tt QUADPACK} code predicts the narrowest accuracy basins in all the
cases. Practically, any computed output $q\sb{2n}$ with actual accuracy
above the computer roundoff is ruled out as unreliable.
At large values of the argument $\omega$ of the trigonometric functions,
this code results in an average rate of spurious outputs of about two percent.
In figs.~\ref{fig:oc10p1} and~\ref{fig:os10p1}, unreliable estimates
of this code are noticed at arguments $\omega > 12\pi$ and $\omega > 14\pi$
respectively.
The user {\sl is not\/} notified of the wrong diagnostics associated to
these outputs at the moment of solving the integrals of interest.
As mentioned in the Introduction, the only way of identifying them is
the far end prediction of nonphysical results for the involved observables.


\begin{table}[ht]\centering
\caption{Comparison of the stability basins and diagnostics reliability
          of the three codes}
\begin{tabular}{|l|cccc|cccc|}
\hline
&\multicolumn{4}{|c|}{Extension $n\sb{max}$}&
 \multicolumn{4}{|c|}{Number of spurious}\\
&\multicolumn{4}{|c|}{of the accuracy basins $\sp{*)}$}&
 \multicolumn{4}{|c|}{diagnostics at output $\sp{**)}$}\\
\hline
Integral family&$\bf (C1)$&$\bf (S1)$&$\bf (C2)$&$\bf (S2)$&
                   $\bf (C1)$&$\bf (S1)$&$\bf (C2)$&$\bf (S2)$\\
\hline
{\tt QUADPACK} &$  219   $&$   192  $&$  359   $&$  379    $&
                   $   254  $&$      33$&$  115   $&$    61  $\\
\hline
Ref.~\cite{AA01}&$  233   $&$   202  $&$  417   $&$  389    $&
                   $     0  $&$      0 $&$    0   $&$     0  $\\
\hline
Present        &$  399   $&$   395  $&$  820   $&$  769    $&
                   $   26   $&$      46$&$   41   $&$     7  $\\
               &$ (521)  $&$  (523) $&$ (1052) $&$ (1069)  $&
                   $   (0)  $&$  (0)   $&$  (0)   $&$    (0) $\\
\hline
\end{tabular}
\end{table}
{\small
$\sp{*)}$ For the present analysis, the upper values correspond to the
          left interrupted line arrows in Figs.~\ref{fig:oc10p1}
	  to~\ref{fig:os11p1}. The values under parentheses correspond
	  to the right interrupted line arrows in the same figures.\\
$\sp{**)}$ For the present analysis, the upper values show the number of
           primary analysis failures. The vanishing values under parentheses
	   show that all the primary analysis failures were corrected
	   under subrange subdivision.
}


The self-validating procedure developed in ref.~%
\cite{AA01}
slightly enlarged the extension of the accuracy basin predictions,
with no wrong outputs at all.

The present reliability analysis identified substantially larger accuracy
basins of the output. All the $q\sb{2n}$ outputs showing more than six
accurate figures have been correctly identified as reliable.
For outputs $q\sb{2n}$ showing inbetween six and three accurate figures,
the present diagnostic was too conservative for 19 $\bf (C1)$ integrals,
51 $\bf (S1)$, 57 $\bf (C2)$, and 176 $\bf (S2)$ integrals.
At values of the argument $\omega$ in large excess to those falling in the
accuracy basins, a number of spurious diagnostics was produced by the primary
reliability analysis. {\sl All\/} the wrong diagnostics occurring at a first
run {\sl were identified as wrong and corrected\/} under subrange subdivision.

\begin{figure}[ht]
\resizebox{\textwidth}{!}
{\includegraphics{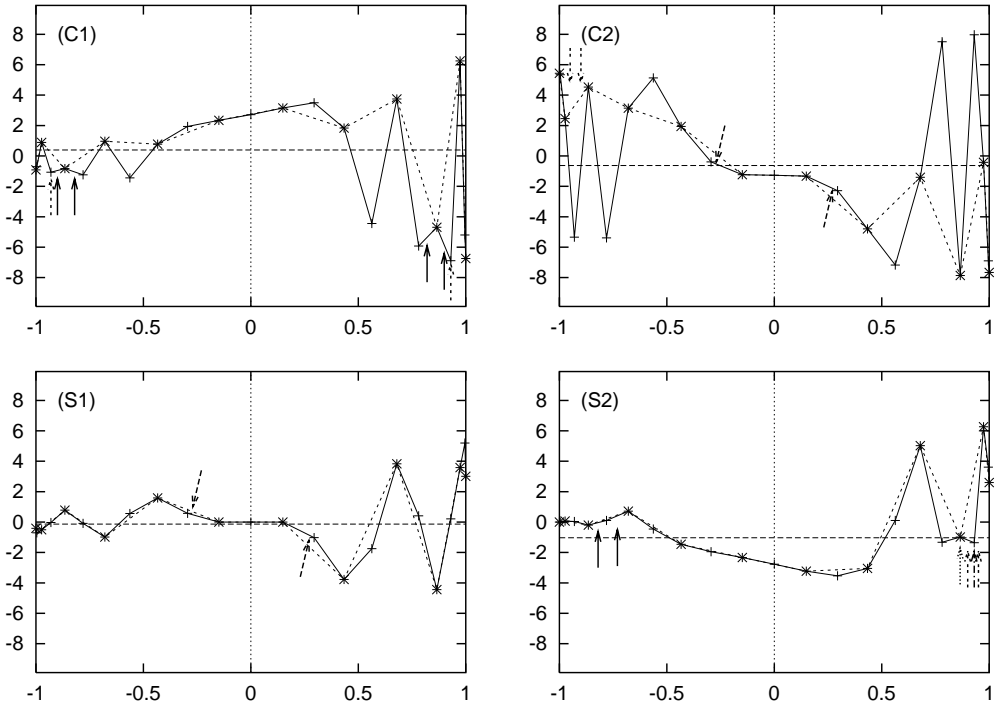}}
\caption{Ill-conditioned integrand features in the family of
         integrals~(\ref{eq:oc10}--\ref{eq:os11}) at $\omega = 1612 \pi / 60$.
        }
\label{fig:1612}
\end{figure}
\begin{figure}[h]
\resizebox{\textwidth}{!}
{\includegraphics{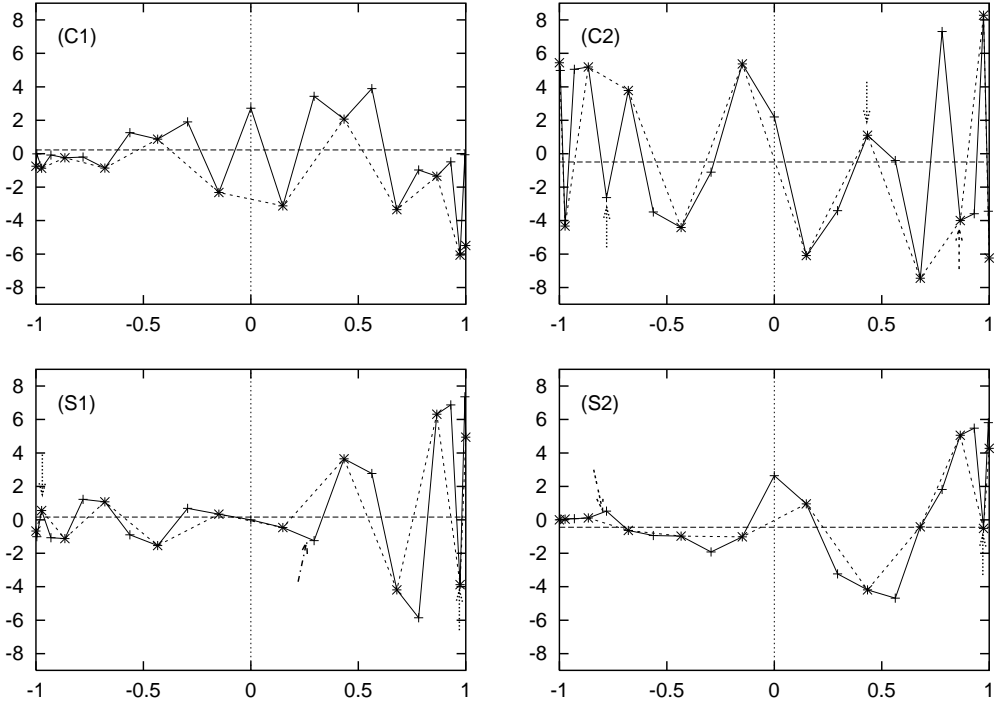}}
\caption{Ill-conditioned integrand features in the family of
         integrals~(\ref{eq:oc10}--\ref{eq:os11}) at $\omega = 3646 \pi / 60$.
        }
\label{fig:3646}
\end{figure}
Fig.~\ref{fig:1612} and fig.~\ref{fig:3646}, show integrand profiles
for the integrals~(\ref{eq:oc10}--\ref{eq:oc11})
and~(\ref{eq:os10}--\ref{eq:os11}) at the large parameter values
$\omega = (1612 \pi / 60)$ and $\omega = (3646 \pi / 60)$ respectively,
together with hints (showed by arrows) on infringements of the reliability
criteria established in Sec.~\ref{sec:wellco}.

A scrutiny of the integrand profiles shows that, in general, it is hardly
probable that a highly oscillatory integrand structure can be resolved at
the existing quadrature knots. However, if intermediate unresolved
structures are present, these induce, as a rule, infringements of one or more
reliability criteria.
The complete list of criteria infringements is given below:
\begin{itemize}
  \item
       {\sl Criterion (I):}\\
       In fig~\ref{fig:1612}:
         the integral $\bf (C1)$ at $y\sb{-8}$ and $y\sb{8}$;
         the integral $\bf (S2)$ at $y\sb{8}$.\\
       In fig~\ref{fig:3646}:
         the integral $\bf (S1)$ at $y\sb{-10}$ and $y\sb{10}$
	   (not shown in the plot);
         the integral $\bf (C2)$ at $y\sb{7}$,
         and the integral $\bf (S2)$ at $y\sb{10}$
           (not shown in the plot).
  \item
       {\sl Criterion (IIa) -- over end subranges:}\\
       In fig~\ref{fig:3646}:
         the integral $\bf (S1)$ over the subrange $[y\sb{-11}, y\sb{-10}]$
	   (not shown in the plot).
  \item
       {\sl Criterion (IIa) -- over inner subranges:}\\
       In fig~\ref{fig:1612}:
         the integral $\bf (C1)$ over the subranges $[y\sb{-8}, y\sb{-7}]$,
	   $[y\sb{-7}, y\sb{-6}]$, $[y\sb{6}, y\sb{7}]$, and
	   $[y\sb{7}, y\sb{8}]$ ;
         the integral $\bf (S2)$ over the subranges $[y\sb{-7}, y\sb{-6}]$
           and $[y\sb{-6}, y\sb{-5}]$;\\
       In fig~\ref{fig:3646}:
         the integral $\bf (C1)$ over 15 (!) subranges: $[y\sb{-10}, y\sb{-9}]$,
	   $[y\sb{-9}, y\sb{-8}]$, $[y\sb{-8}, y\sb{-7}]$,
	   $[y\sb{-7}, y\sb{-6}]$, $[y\sb{-6}, y\sb{-5}]$,
	   $[y\sb{-4}, y\sb{-3}]$, $[y\sb{-3}, y\sb{-2}]$,
	   $[y\sb{-1}, y\sb{0}]$, $[y\sb{2}, y\sb{3}]$; $[y\sb{3}, y\sb{4}]$,
	   $[y\sb{5}, y\sb{6}]$, $[y\sb{6}, y\sb{7}]$, $[y\sb{7}, y\sb{8}]$,
	   $[y\sb{8}, y\sb{9}]$, and $[y\sb{9}, y\sb{10}]$
	   (not shown in the plot since they are obvious).
  \item
       {\sl Criterion (IIb):}\\
       In fig~\ref{fig:1612}:
         the integral $\bf (S2)$ at $y\sb{7}$;\\
       In fig~\ref{fig:3646}:
         the integral $\bf (S1)$ at $y\sb{-9}$ and $y\sb{9}$;
         the integral $\bf (C2)$ at $y\sb{-6}$ and $y\sb{3}$;
         the integral $\bf (S2)$ at $y\sb{9}$.
  \item
       {\sl Criterion (III):}\\
       In fig~\ref{fig:1612}:
         the integral $\bf (C2)$ inside the subranges: $[y\sb{-9}, y\sb{-8}]$,
	   and $[y\sb{-8}, y\sb{-7}]$;
         the integral $\bf (S2)$ inside the subranges: $[y\sb{7}, y\sb{8}]$,
	   and $[y\sb{8}, y\sb{9}]$;
  \item
       {\sl Criterion (IV):}\\
       In fig~\ref{fig:3646}:
         the integral $\bf (S1)$ at $y\sb{9}$ (left derivative).
  \item
       {\sl Criterion (V):}\\
       In fig~\ref{fig:3646}:
         the integral $\bf (S2)$ at $y\sb{-6}$ (left neighbourhood).
  \item
       {\sl Criterion (VI):}\\
       In fig~\ref{fig:1612}:
         the integral $\bf (S1)$ to the right of the knot $y\sb{-2}$ and
	   to the left of the knot $y\sb{2}$;
         the integral $\bf (C2)$ to the right of the knot $y\sb{-2}$ and
	   to the left of the knot $y\sb{2}$.
\end{itemize}
\section{Comments and conclusions}
\label{sec:concl}

The present investigation started from the need to get reliable numerical
solutions of difficult parametric integrals occurring in theoretical models
devoted to the study of the mechanism of the high-$T\sb{c}$ superconductivity
in cuprates
\cite{Plak95}-\cite{YuLu00}
and in a theoretical model of nuclear fission
\cite{Aurel98}.
An important prerequisite to be satisfied by the automatic quadrature
algorithms needed for the evaluation of the observables was the substantial
increase of the reliability of the local error estimates.

We have found that the study of the conditioning of the integrand profile
enables the formulation of validation criteria (consistency conditions for
a well-conditioned profile) able to identify insufficient profile resolution
or the occurrence of isolated difficult points of the integrand.
The analysis is simple, it is intuitive, it is easily implemented in a
computer program and it is easily done.

An important supplementary bonus offered by this analysis was the
identification of $q\sb{2n}$ output reliability ranges which are
substantially larger in comparison with those obtained within the 
usual implementations of quadrature routines.
 The unsatisfactory features noticed in the validation criteria developed in
ref.~%
\cite{AA01}
have been fully removed.

The subroutines doing the profile analysis described in this paper are
documented and described in a separate document
\cite{AA02}.

We conclude this study with the observation that
the validation analysis described in the present paper is {\sl not\/}
intended to replace the existing quadrature algorithms. When the estimated
accuracy exceeds a critical threshold (tentatively set to five decimal
figures), then the present procedure is skipped altogether.
  However, if this threshold is not attained, it is automatically activated
by the general control routine. Its results prove to be invaluable in the
analysis of complex integrands, where it is able to discover the
overwhelming fraction of peculiar integrand profiles at early stages of
the analysis.
\section*{Acknowledgments}
     The investigation was partially financed by
     the JINR grant no.~571/15.10.2001 afforded by the Romanian
     Plenipotentiary Representative.

     One author (Gh.A.) is grateful to Yu Lu, A. S\v andulescu, and
     \c S. Mi\c sicu for discussions of specific physical models.


\vfil\eject
\def\endpage{\hfill\eject}

\end{document}